\documentclass[12pt,a4paper]{article}

\usepackage{graphicx}
\usepackage{color}
\usepackage{amsmath}
\usepackage{amssymb}
\usepackage{amscd}
\usepackage{amsthm}
\usepackage{delarray}
\usepackage{enumerate}
\usepackage[T1]{fontenc}
\usepackage{inputenc}
\usepackage{enumerate}
\input xy

\begin{document}

\xyoption{all}
\setlength{\parindent}{5mm}
\renewcommand{\leq}{\leqslant}
\renewcommand{\geq}{\geqslant}
\newcommand{\N}{\mathbb{N}}
\newcommand{\sph}{\mathbb{S}}
\newcommand{\Z}{\mathbb{Z}}
\newcommand{\R}{\mathbb{R}}
\newcommand{\C}{\mathbb{C}}
\newcommand{\F}{\mathbb{F}}
\newcommand{\g}{\mathfrak{g}}
\newcommand{\h}{\mathfrak{h}}
\newcommand{\K}{\mathbb{K}}
\newcommand{\RN}{\mathbb{R}^{2n}}
\newcommand{\derive}[2]{\frac{\partial{#1}}{\partial{#2}}}
\renewcommand{\S}{\mathbb{S}}
\renewcommand{\H}{\mathbb{H}}
\newcommand{\eps}{\varepsilon}
\theoremstyle{plain}
\newtheorem{theo}{Theorem}
\newtheorem{prop}[theo]{Proposition}
\newtheorem{lem}[theo]{Lemma}
\newtheorem{definit}[theo]{Definition}
\newtheorem{corol}[theo]{Corollary}
\newcommand{\preuve}{\textit{Proof: }}

\title{Hamiltonian pseudo-representations}
\author{V. Humili\`ere}
\date{20/07/07}
\maketitle
\normalsize
\begin{center}
Centre de Mathématiques Laurent Schwartz

UMR 7640 du CNRS

Ecole Polytechnique - 91128 Palaiseau, France

\texttt{vincent.humiliere@math.polytechnique.fr}
\end{center}

\begin{abstract}The question studied here is the behavior of the Poisson
bracket under $C^0$-perturbations. In this purpose, we introduce the
notion of pseudo-representation and prove that the limit of a
converging pseudo-representation of any normed Lie algebra is a
representation.

An unexpected consequence of this result is that for many non-closed
symplectic manifolds (including cotangent bundles), the group of
Hamiltonian diffeomorphisms (with no assumptions on supports) has no
$C^{-1}$ bi-invariant metric. Our methods also provide a new proof
of Gromov-Eliashberg Theorem, it is to say that the group of
symplectic diffeomorphisms is $C^0$-closed in the group of all
diffeomorphisms.
\end{abstract}

\section{Statement of results}

\subsection{Poisson Brackets and $C^0$-convergence}

We consider a symplectic manifold $(M,\omega)$. A function $H$ on
$M$ will be said normalized if $\int_M{H}\omega^n=0$ for $M$ closed
or if $H$ has compact support otherwise. We will denote
$C_0^{\infty}(M)$ the set of normalized smooth functions. Endowed
with the Poisson brackets $\{\cdot,\cdot\}$, it has the structure of
a Lie algebra.

In the whole paper, we will denote $X_H$ the symplectic gradient of
a smooth function $H$, i.e., the only vector field satisfying
$dH=\iota_{X_H}\omega$. Then, the Poisson brackets are given by
$\{H,K\}=dH(X_K)$.

Let $\g$ be a normed Lie algebra, i.e., a Lie algebra endowed with a
norm $\|\cdot\|$ such that for some constant $C$,
$$\|[f,g]\|\leq C\|f\|\cdot\|g\|,$$ and consider the following
definition.

\begin{definit}A sequence of linear maps
$$\rho_n:(\g,\|\cdot\|)\,\to\,(C_0^{\infty}(M),\|\cdot\|_{C^0}),$$
will be called a \textit{pseudo-representation} if the sequence of
bilinear maps
$$B_n:(f,g)\mapsto\{\rho_n(f),\rho_n(g)\}-\rho_n([f,g])$$
converges to 0.
\end{definit}

If it has a limit, we may ask whether this limit is a
representation. If so, we would have
$$\{\rho_n(f),\rho_n(g)\}\to\{\rho(f),\rho(g)\},\text{ for all
}f,g\in\g.$$ This has been proved in \cite{CV} for abelian Lie
algebras. The main result of this paper is that it holds for all
normed Lie algebras.

\begin{theo}\label{theo principal} For any normed Lie algebra (in particular for finite dimensional Lie
algebras), the limit of a converging pseudo-representation is a
representation.
\end{theo}

\medskip
\noindent \textit{Remark 1.} This result generalizes
Gromov-Eliashberg's Theorem of $C^0$ closure of the
symplectomorphisms group in the group of diffeomorphisms.

Indeed, a diffeomorphism of $\RN$ is symplectic if and only if its
coordinate functions $(f_i), (g_i)$ satisfy
$$\{f_i,g_j\}=\delta_{ij},\quad \{f_i,f_j\}=\{g_i,g_j\}=0.$$
Thus we can easily see that a sequence of symplectomorphisms gives a
pseudo-representation of a 2-nilpotent Lie algebra. If the support
of the coordinate functions were compact, we could immediately apply
Theorem \ref{theo principal}. In fact, for compactly supported
symplectomorphisms, these functions are affine at infinity, and we
have to adapt the proof to this case (See Appendix \ref{groelia} for
details).

\medskip \noindent \textit{Remark 2.} Consider the following
question: If $F_n$, $G_n$ and $\{F_n,G_n\}$ respectively converge to
$F$, $G$ and $H$ (all function being smooth and normalized, and all
convergence being in the $C^0$ sense), is it true that $\{F,G\}=H$ ?

Theorem \ref{theo principal} states that the answer is positive when
there is some Lie algebra structure. Nevertheless, in general, the
answer is negative, as shows the following example, which is derived
from Polterovich's example presented in Section \ref{preuve 2}. Let
$\chi$ be a compactly supported smooth function on $\R$, and set the
following functions on $\R^2$:
$$F_n(q,p)=\frac{\chi(p)}{\sqrt{n}}\cos(nq),$$
$$G_n(q,p)=\frac{\chi(p)}{\sqrt{n}}\sin(nq).$$
It is easy to see that $F_n$ and $G_n$ converge to 0, but that their
Poisson brackets equal $\chi(p)\chi'(p)\neq 0$.

This example shows that when the Poisson brackets $C^0$-converge,
then its limit is not necessarily the brackets of the respective
limits. But in that case, we can see that the Hamiltonians $F_n$ and
$G_n$ do not generate a pseudo-representation.

\medskip \noindent\textit{Remark 3.} The theorem holds if we replace
the symplectic manifold with a general Poisson manifold. Indeed,
Poisson manifolds are foliated by Poisson submanifolds that are
symplectic, and we just have to apply theorem \ref{theo principal}
to each leaf.

\medskip \noindent\textit{Remark 4.} The theorem leads us to the
following \begin{definit} A \textnormal{continuous Hamiltonian
representation} of a normed Lie algebra $\mathfrak{g}$ is a
continuous linear map $\mathfrak{g}\to C^0(M)$ which is the
$C^0$-limit of some pseudo-representation of $\mathfrak{g}$.
\end{definit}
We will not study this notion further in this paper. Nevertheless
let us give some example:

\noindent\textit{Example: }Let $\rho:\g\to C_0^\infty(M)$ be a
smooth Hamiltonian representation in the usual sense, and let
$\varphi$ be a homeomorphism of $M$ which is the $C^0$-limit of a
sequence of symplectomorphisms. Then, $\rho':\g\to C^0(M)$, given by
$\rho'(g)=\rho(g)\circ\varphi$, is clearly a continuous Hamiltonian
representation.

\medskip
\noindent \textit{Question 1:} Given two sequences of Hamiltonians
$(F_n)$, $(G_n)$ that $C^0$-converge to smooth $F$ and $G$, is there
some sufficient condition for the bracket $\{F,G\}$ not to be the
limit of the brackets $\{F_n,G_n\}$? Propositions
\ref{distributions1} and \ref{distributions2} give restrictions on
the possible counter- examples.

\medskip
\noindent \textit{Question 2:} Let us consider the following number
introduced by Entov, Polterovich and Zapolsky in \cite{EPZ}:
$$\Upsilon(F,G)=\liminf_{\eps\to 0}\left\{\|\{F',G'\}\|\,|\, \|F-F'\|_{C^0}<\eps,\|G-G'\|_{C^0}<\eps\right\}$$
The result of Cardin and Viterbo mentioned above which is exactly
Theorem \ref{theo principal} in the abelian case can be restated as
follows:
$$\Upsilon(F,G)>0 \text{ if and only if }\{F,G\}\neq 0.$$
Entov, Polterovich and Zapolsky have improved this result by giving
explicit lower bounds on $\Upsilon(F,G)$, in terms of quasi-states
(see \cite{EPZ} and \cite{Z}). We may wonder whether there exist
similar inequalities in the non abelian case.

\subsection{Bi-invariant Metrics}

Here we consider a subgroup $\mathcal{G}$ of the group
$\mathcal{H}(M)$ of Hamiltonian diffeomorphisms on $M$. If we denote
$\phi_H^t$ the flow generated by $X_H$ (when it exists), and
$\phi_H=\phi_H^1$ the time-1 map, $\mathcal{H}(M)$ is the set of all
diffeomorphisms $\phi$ for which it exists a path of Hamiltonian
functions $H_t\in C^{\infty}(M)$ such that $\phi=\phi_H$.

\begin{definit} A \textnormal{bi-invariant metric} on $\mathcal{G}$ is a distance $d$ on $\mathcal{G}$ such that for any $\phi,\psi,\chi$ in $\mathcal{G}$,
$$d(\phi,\psi)=d(\phi\chi,\psi\chi)=d(\chi\phi,\chi\psi).$$
It will be said $C^{-1}$ if its composition with the map
$\Phi:H\mapsto\phi_H^1$ is a continuous map
$\Phi^{-1}(\mathcal{G})\times\Phi^{-1}(\mathcal{G})\to\R$, where
$\Phi^{-1}(\mathcal{G})\subset Ham$ is endowed with the compact-open
topology.
\end{definit}

There are several well known examples of $C^{-1}$ bi-invariant
metrics, as, for example, Hofer's metric defined on the subgroup
Hamiltonian diffeomorphisms generated by compactly supported
functions $\mathcal{H}_c(M)$ (see \cite{H1} or \cite{HZ}), Viterbo's
metric defined on $\mathcal{H}_c(\RN)$ (see \cite{V1}), and its
analogous version defined by Schwarz in \cite{Sc} for symplectically
aspherical closed symplectic manifolds.

As far as we know, if we remove the assumption of compactness of the
support, the question whether there exists such metrics is still
open. Here we prove that the answer is negative for a large class of
symplectic manifolds.

Let $(N,\xi)$ be a contact manifold with contact form $\alpha$
(i.e., a smooth manifold $N$ with a smooth hyperplane section $\xi$
which is locally the kernel of a 1-form $\alpha$ whose differential
$d\alpha$ is non-degenerate on $\xi$). Its \textit{symplectization}
is by definition the symplectic manifold $\mathcal{S}N=\R\times N$
endowed with the symplectic form $\omega=d(e^s\alpha)$, where $s$
denotes the $\R$-coordinate in $\R\times N$. For any contact form
$\alpha$, one can define the \textit{Reeb vector field} $X_R$ by the
identities $\iota_{X_R} d\alpha$ and $\alpha(X_R)=1$. The
trajectories of $X_R$ are called \textit{characteristics}. The
question of the existence of a closed characteristic constitutes the
famous Weinstein's conjecture. It has now been proved for large
classes of contact manifolds (see e.g.
\cite{FHV,HV1,HV2,Lu,LT,V2,T}...).

Let us now state our result that will be proved in section
\ref{preuve 2}

\begin{theo}\label{theo principal 2}If $M$ is the symplectization of a contact manifold whose dimension is at least $3$ and that admits a closed characteristic, then there is no $C^{-1}$ bi-invariant metric on
$\mathcal{H}(M)$.
\end{theo}

\begin{corol}\label{cotangents}If $N$ is a smooth manifold whose dimension is at least 2 and if $T^*N$ is its cotangent bundle, then there is no $C^{-1}$ bi-invariant metric on
$\mathcal{H}(T^*N)$.
\end{corol}

\medskip
\noindent \textit{Remark.} At least in the case of manifolds of
finite volume, there probably exists non closed manifolds with such
distances. Indeed, it follows from our previous work \cite{Hu} that
 Viterbo's metric extends to Hamiltonians functions smooth out of a
 "small" compact set. Replacing Viterbo's metric with Schwarz's metric, we can
 reasonably expect to have: If $M^{2n}$ is a closed symplectically
 aspherical manifold and $K$ is a closed submanifold of dimension $\leq
 n-2$, then Schwarz's metric on $\mathcal{H}(M)$ extends to
 $\mathcal{H}(M-K)$.

\section{Proofs}
\subsection{Identities for Hamiltonian pseudo-representations}

\medskip
\begin{lem}\label{limite} Let $\rho_n$ be a bounded (not necessarily converging) pseudo-representation of a normed Lie algebra $\g$. Let $f,g\in\g$, then the sequence of Hamiltonian functions
$$\rho_n(f)\circ\phi_{\rho_n(g)}^{s}\,-\,\sum_{j=0}^{+\infty}\rho_n(ad(g)^jf)\,\frac{s^j}{j!}$$
converges to zero for the $C^0$-norm on $M$. Moreover, the
convergence is uniform over the $s$'s in any compact interval.
\end{lem}

\medskip
\noindent \textit{Remark:} For a representation equality holds. It
recalls the Baker-Campbell-Haussdorf formula.

\medskip
\noindent \preuve First remark that the considered sum converges.
Indeed, the $C^0$-norm of its remainder can be bounded by the
remainder of a converging sum, as follows:
$$\left\|\sum_{j=N}^{+\infty}\rho_n(ad(g)^jf)\,\frac{s^j}{j!}\right\|\leq \sum_{j=N}^{+\infty}R\|f\|\frac{(sC\|g\|)^j}{j!}.$$
where $R$ is an $n$-independent upper bound for the sequence
$$\|\rho_n\|=\sup\{\|\rho_n(h)\|_{C^0}\,|\,\|h\|=1\}.$$

Now, let us prove our lemma. Poisson equation gives
$$\frac{d}{ds}(\rho_n(f)\circ\phi_{\rho_n(g)}^{s})=\{\rho_n(f),\rho_n(g)\}\circ\phi_{\rho_n(g)}^{s}$$
and hence
$$
\rho_n(f)\circ\phi_{\rho_n(g)}^{s_0}= \rho_n(f) +
\int_0^{s_0}\{\rho_n(f),\rho_n(g)\}\circ\phi_{\rho_n(g)}^{s_1}\,ds_1$$

$$=\rho_n(f)+\int_0^{s_0}\rho_n([f,g])\circ\phi_{\rho_n(g)}^{s_1}\,ds_1+ \int_0^{s_0}B_n(f,g)\circ\phi_{\rho_n(g)}^{s_1}\,ds_1.$$

Then, by a simple induction, we get for all integer $N$:
\begin{equation*}\label{prem}
\rho_n(f)\circ\phi_{\rho_n(g)}^{s_0}=\sum_{j=0}^N\rho_n(ad(g)^jf)\,\frac{{s_0}^j}{j!}\,+\,R_{N,n}(s_0)\,+\,S_{N,n}(s_0),
\end{equation*}
where,
$$R_{N,n}(s_0)=\int_0^{s_0}\int_0^{s_1}\cdots\int_0^{s_{N}}\rho_n(ad(g)^{N+1}f)\circ\phi_{\rho_n(g)}^{s_{N+1}}ds_{N+1}\cdots ds_1$$
\begin{eqnarray*}
S_{N,n}(s_0)=\sum_{j=0}^{N}
\int_0^{s_0}\int_0^{s_1}\cdots\int_0^{s_{j}}B_n(ad(g)^{j}f,g))\circ\phi_{\rho_n(g)}^{s_{j+1}}ds_{j+1}\cdots
ds_1
 \end{eqnarray*}

\medskip

Let us now denote
$$\|B_n\|=\sup\{\|\{\rho_n(f),\rho_n(g)\}-\rho_n([f,g])\|_{C^0}\,|\,\|f\|=\|g\|=1\}.$$
By assumptions $\|B_n\|$ converges to 0.

\medskip

Then,
\begin{eqnarray*}
\|R_{N,n}(s_0)\|_{C^0} & \leq &
\int_0^{s_0}\int_0^{s_1}\cdots\int_0^{s_{N-1}}R\|g\|^NC^N\|f\|ds_N\cdots
ds_1,\\
& \leq & R\|f\|\frac{\|g\|^NC^Ns_0^N}{N!}, \end{eqnarray*} which
proves that $R_{N,n}(s_0)$ converges to $0$ with $N$, uniformly in
$n$.

In addition,
\begin{eqnarray*}
\|S_{N,n}(s_0)\|\leq \sum_{j=0}^{N-2}
\int_0^{s_0}\int_0^{s_1}\cdots\int_0^{s_{j}}\|B_n\|\|f\|\|g\|^{j}ds_{j+1}\cdots
ds_1
 \end{eqnarray*}
We thus have $\|S_{N,n}(s_0)\|\leq \|B_n\|\,\|f\|\exp(s_o\|g\|)$ for
any $N$. As a consequence, letting $N$ converge to $+\infty$, we get
$$\left\|\rho_n(f)\circ\phi_{\rho_n(g)}^{s}\,-\,\sum_{j=0}^{+\infty}\rho_n(ad(g)^jf)\,\frac{s^j}{j!}\right\|\leq\|B_n\|\,\|f\|\exp(s_o\|g\|).$$
This achieves the proof because the right hand side converges to
$0$.$\quad\Box$

\subsection{Proof of theorem \ref{theo principal}}

Let $f,g\in\g$. We want to prove that
$\{\rho(f),\rho(g)\}=\rho([f,g])$. We can assume without loss of
generality that $\|g\|<1$.

By Lemma \ref{limite},
$$\rho_n(f)\circ\phi_{\rho_n(g)}^{s}\,-\,\sum_{j=0}^{+\infty}\rho_n(ad(g)^jf)\,\frac{s^j}{j!}\stackrel{C^0}{\to}0.$$
Each term of the sum converges with $n$. Since the sum converges
uniformly in $n$, we get that for any $s$,
$$\rho_n(f)\circ\phi_{\rho_n(g)}^{s}\,\stackrel{C^0}{\to}\,\sum_{j=0}^{+\infty}\rho(ad(g)^jf)\,\frac{s^j}{j!}.$$
As a consequence, the flow generated by
$\rho_n(f)\circ\phi_{\rho_n(g)}^{s}$ $\gamma$-converges to the flow
generated by $\sum_{j=0}^{+\infty}\rho(ad(g)^jf)\,\frac{s^j}{j!}$.

But on the other hand, the flow of
$\rho_n(f)\circ\phi_{\rho_n(g)}^{s}$ is
$t\mapsto\phi_{\rho_n(g)}^{-s}\phi_{\rho_n(f)}^{t}\phi_{\rho_n(g)}^{s}$,
which $\gamma$-converges to
$\phi_{\rho(g)}^{-s}\phi_{\rho(f)}^{t}\phi_{\rho(g)}^{s}$. Indeed,
$\rho_n(g)\stackrel{C^0}{\to}\rho(g)$ and
$\rho_n(f)\stackrel{C^0}{\to}\rho(f)$ which implies that there
respective flow $\gamma$-converges.

Therefore,
$t\mapsto\phi_{\rho(g)}^{-s}\phi_{\rho(f)}^{t}\phi_{\rho(g)}^{s}$ is
the flow of $\sum_{j=0}^{+\infty}\rho(ad(g)^jf)\,\frac{s^j}{j!}$.
The functions being normalized,
$$\rho(f)\circ\phi_{\rho(g)}^{s}=\sum_{j=0}^{+\infty}\rho(ad(g)^jf)\,\frac{s^j}{j!}.$$
Now, first taking derivative with respect to $s$, we get
$\{\rho(f),\rho(g)\}=\rho([f,g])$. $\quad\Box$

\subsection{Proof of theorem \ref{theo principal 2}}\label{preuve 2}

Let us consider the following Hamiltonian functions on $\R^2$ (this
example is due to Polterovich) with symplectic form written in polar
coordinates $rdr\wedge d\theta$.
$$F_n(r,\theta)=\frac{r}{\sqrt{n}}\cos(n\theta),$$
$$G_n(r,\theta)=\frac{r}{\sqrt{n}}\sin(n\theta).$$
We see that $\{F_n,G_n\}=1$ and that $F_n$ and $G_n$ converge to
$0$. Now, consider $\g$ the 3-dimensional Heisenberg Lie algebra
(i.e., the Lie algebra with basis $\{f,g,h\}$ such that $[f,g]=h$
and $[f,h]=[g,h]=0$) and set $\rho_n(f)=F_n$, $\rho_n(g)=G_n$ and
$\rho_n(h)=1$. Then, $\rho_n$ is a pseudo-representation of $\g$ in
$Ham(R^2)$. The limit $\rho$ of $\rho_n$ satisfies $\rho(f)=0$,
$\rho(g)=0$, $\rho(h)=1$. Since $\{\rho(f),\rho(g)\}\neq\rho(h)$,
$\rho$ is not a representation of $\g$.

Since $\g$ has finite dimension, this example shows that Theorem
\ref{theo principal} is false in general if we replace
$C_0^{\infty}(M)$ with $C^{infty}(M)$ for a non-compact manifold
$M$, and uniform convergence with the uniform convergence on compact
sets (compact-open topology).

If we read carefully the proof of Theorem \ref{theo principal}, we
see that the whole proof can be repeated in this settings except the
three following points where the compactness of supports are needed
\begin{itemize}
  \item Each time we consider the flows of the Hamiltonians, they must
  be complete. This is automatic for compactly supported Hamiltonians, but false in
  general. With the notations of the proof, the flows needed are
  those of $\rho_n(f)$, $\rho(f)$, $\rho_n(g)$, $\rho(g)$ and $\sum_{j=0}^{+\infty}\rho(ad(g)^jf)\,\frac{s^j}{j!}$.
  \item The functions $\rho_n(f)$, $\rho(f)$, $\rho_n(g)$, $\rho(g)$
  have to be normalized in some sense.
  \item We use a $C^{-1}$ bi-invariant metric. This exists on $\mathcal{H}_c(M)$, but we do not know whether it exists on $\mathcal{H}(M)$.
\end{itemize}

The following lemma follows from the above discussion.

\begin{lem}\label{précis}Let $M$ be a non-compact symplectic manifold, $\g$ a normed Lie algebra, and $\rho_n$ a pseudo-representation of $\g$ in
$Ham(M)$, with limit $\rho$. Suppose there exists two elements $f$
and $g$ in $\g$, such that:
\begin{itemize}
\item all the Hamiltonian functions $\rho_n(f)$, $\rho(f)$, $\rho_n(g)$, $\rho(g)$ and $\sum_{j=0}^{+\infty}\rho(ad(g)^jf)\,\frac{s^j}{j!}$ exist and have complete flows,
\item there exists an open set on which all the functions $\rho_n(f)$, $\rho(f)$, $\rho_n(g)$,
$\rho(g)$ vanish identically.
\item $\{\rho(f),\rho(g)\}\neq\rho([f,g])$.
\end{itemize}
Then the group of Hamiltonian diffeomorphisms $\mathcal{H}(M)$
admits no $C^{-1}$ bi-invariant metric.$\quad\Box$
\end{lem}

\medskip
\noindent \textit{Proof of Theorem \ref{theo principal 2}: }We want
to apply Lemma \ref{précis}. We first consider the case of $\S^1$.
In that case we are not able to get the second requirement of Lemma
\ref{précis}, but let us show how we get the others.

We just adapt Polterovich's example by setting :
$$\rho_n(f)(s,\theta)=\frac{e^{s/2}}{\sqrt{n}}\cos(n\theta),$$
$$\rho_n(g)(s,\theta)=\frac{e^{s/2}}{\sqrt{n}}\sin(n\theta).$$
The symplectic form being defined on $\R\times \S^1$ by
$d(e^sd\theta)=e^sds\wedge d\theta$, we get
$\{\rho_n(f),\rho_n(g)\}=2$. Since $\rho(f)=\rho(g)=0$ we have a
pseudo-representation of the 3-dimensional Heisenberg Lie algebra,
and its limit is not a representation. We can also verify that all
elements $\rho_n(f)$, $\rho(f)$, $\rho_n(g)$, $\rho(g)$ and
$\sum_{j=0}^{+\infty}\rho(ad(g)^jf)\,\frac{s^j}{j!}$ exist and have
complete flows for $f$, $g$ generators of the 3-dimensional
Heisenberg Lie algebra, and $\rho_n,\rho$ as in the example.

Since $\rho(f)=0$, $\rho(g)=0$ and
$\sum_{j=0}^{+\infty}\rho(ad(g)^jf)\,\frac{s^j}{j!}=2s$, this is
obvious for them.

The Hamiltonian vector field of $\rho_n(f)$ is
$$\left(e^{-s/2}\sqrt{n}\sin(n\theta)\right)\derive{}{\theta}-
\left(\frac{1}{2\sqrt{n}}e^{-s/2}\cos(n\theta)\right)\derive{}{s},$$
which is equivalent through the symplectomorphism
$$(\R\times \S^1,d(e^sd\theta))\to(\R^2-\{0\}, rdr\wedge d\theta)),\
(s,\theta)\mapsto(e^{-s/2},\theta)),$$ to the vector field
$$\left(r\sqrt{n}\sin(n\theta)\right)\derive{}{\theta}+
\left(\frac{1}{\sqrt{n}}\cos(n\theta)\right)\derive{}{r}.$$ The norm
of this vector field is bounded by a linear function in $r$.
Therefore, it is a consequence of Gronwall's lemma that it is
complete.

\medskip

Let us consider now the case $d=\dim(N)\geq 3$. There, we will be
able to get all the requirements of Lemma \ref{précis}. Denote by
$\gamma$ a closed characteristic, parameterized by $\theta\in\S^1$.
Since the Reeb vector field is transverse to the contact structure
$\xi$, there exists a diffeomorphism that maps a neighborhood
$\mathcal{V}_0$ of the zero section in the restricted bundle
$\xi|_{\gamma}$, onto a neighborhood $\mathcal{V}_1$ of $\gamma$ in
the contact manifold $N$. Since $\xi|_{\gamma}$ is a symplectic
bundle over $\S^1$, it is trivial. We thus have a neighborhood $U$
of $0$ in $\RN$ and a diffeomorphism $\psi:\S^1\times U \to
\mathcal{V}_1\subset N$. The pull back of $\xi$ by $\psi$ is a
contact structure on $\S^1\times U$ which is contactomorphic (via
Moser's argument) to the standard contact structure $d\theta-pdq$ on
$\S^1\times U$. Therefore, the above diffeomorphism $\psi$ can be
chosen as a contactomorphism.


Then the symplectization $\mathcal{S}\gamma$ of the closed
characteristic gives a symplectic embedding
$\mathcal{S}\S^1\hookrightarrow\mathcal{S}N$. This embedding admits
$\mathcal{S}(\S^1\times U)$ as a neighborhood. Moreover, if we
denote $s$, $\theta$ and $x$ the coordinates in
$\mathcal{S}(\S^1\times U)$, $\psi$ has been constructed so that $s$
and $\theta$ are conjugated variables and the direction of $x$ is
symplectically orthogonal to those of $s$ and $\theta$. That will
allow the following computations.

Just like in the above example, we have a pseudo-representation of
$\g$ if we consider
$$(\rho_n(f))(s,\theta,x)=\frac{\chi(x)e^{s/2}}{\sqrt{n}}\cos(n\theta),$$
\begin{equation}\label{lexemple}(\rho_n(g))(s,\theta,x)=\frac{\chi(x)e^{s/2}}{\sqrt{n}}\sin(n\theta),\end{equation}
and $(\rho_n(h))(s,\theta,x)=2\chi(x)^2$. Indeed, we have again
$\{\rho_n(f),\rho_n(g)\}=\rho_n(h)$, but its limit $\rho$ satisfies
$\{\rho(f),\rho(g)\}=0\neq 1=\rho(h)$ and is not a representation.
The fact that the elements $\rho_n(f)$, $\rho(f)$, $\rho_n(g)$,
$\rho(g)$ and $\sum_{j=0}^{+\infty}\rho(ad(g)^jf)\,\frac{s^j}{j!}$
exist and have complete flows follows from the case
$d=1$.$\quad\Box$

\medskip
\noindent\textit{Proof of Corollary \ref{cotangents}} Let $M$ be a
smooth manifold, and choose a Riemannian metric on it. Then,
consider the symplectization $\mathcal{S}\S T^*M$ of the sphere
cotangent bundle $\S T^*M$. The cotangent bundle can be seen as the
compactification of $\mathcal{S}\S T^*M$, the set at infinity being
the zero section of $T^*M$ (or $\{-\infty\}\times\S T^*M$ if we see
$\mathcal{S}\S T^*M$ as $\R\times\S T^*M$).

The Reeb flow of $\S T^*M$ projects itself to the geodesic flow on
$M$, and the closed characteristics are exactly the trajectories
that project themselves to closed geodesics. Since any closed
manifold carries a closed geodesic(see \cite{Kl}), we can consider
Example (\ref{lexemple}). It clearly extends to the compactification
(the Hamiltonian functions involved and all their derivatives
converges to $0$ when $s$ goes to $-\infty$), and we can achieve the
proof as for Theorem \ref{theo principal 2}.$\quad\Box$

\appendix

\section{A proof of Gromov-Eliashberg theorem.}\label{groelia}

In this section, we show how our methods allow to recover
Gromov-Eliashberg Theorem.

\begin{theo}[Gromov, Eliashberg] The group of compactly supported symplectomorphisms
$Symp_c(\RN)$ is $C^0$-closed in the group of all diffeomorphisms of
$\RN$.
\end{theo}

\noindent\textit{Proof.} Let $\phi_n$ be a sequence of
diffeomorphisms that converges uniformly to a diffeomorphism $\phi$.
Denote $(f^n_i), (g^n_i)$ (resp. $f_i, g_i$) the coordinate
functions of $\phi_n$ (resp. $\phi$). These coordinate functions can
be seen has Hamiltonian functions affine at infinity (i.e., that can
be written $H+u$ with $H\in Ham_c$ and $u$ affine map). Moreover,
for a given sequence $(f^n_i)$ or $(g^n_i)$, the linear part does
not depend on $n$.

Since $\phi_n$ is symplectic, we have:
$$\{f^n_i,g^n_j\}=\delta_{ij},\quad \{f^n_i,f^n_j\}=\{g^n_i,g^n_j\}=0.$$
Thus the coordinate functions of $\phi_n$ give a
pseudo-representation of the 2-nilpotent Lie algebra $\g$ generated
by elements $a_i, b_i, c$, with the relations
$$[a_i,b_j]=\delta_{ij},\ [a_i,a_j]=[b_i,b_j]=0, \text{and}\ [a_i,c]=[b_i,c]=0.$$
Since $\phi$ is symplectic if and only if
$$\{f_i,g_j\}=\delta_{ij},\quad \{f_i,f_j\}=\{g_i,g_j\}=0$$
the proof will be achieved if we prove that the limit of this
pseudo-representation is a representation. Consequently, we have to
adapt the proof of Theorem \ref{theo principal} to the case of
Hamiltonian functions affine at infinity, for 2-nilpotent Lie
algebras. Gromov-Eliashberg Theorem then follows from the next two
lemmas.

\begin{lem}\label{lineaire1} Let $u$, $v$ be two affine maps
$\RN\to\R$ and $H_n$, $K_n$ be compactly supported Hamiltonians,
such that $$H_n\to H, \ K_n\to K,\ \{H_n+u,K_n+v\}\to 0.$$ Then
$\{H+u,K+v\}=0$.
\end{lem}

\begin{lem}\label{lineaire2}Let $u$, $v$, $w$ be linear forms on $\RN$, and $H_n$, $K_n$, $G_n$, be compactly supported Hamiltonians such that
$$H_n\to H,\ K_n\to K,\ G_n\to G,$$
$$\{H_n+u,G_n+w\}\to 0,$$
$$\{K_n+v,G_n+w\}\to 0,$$
$$\{H_n+u,K_n+v\}-(G_n+w)\to 0.$$
Then $\{H+u,G+w\}=0$, $\{K+v,G+w\}=0$ and $\{H+u,K+v\}=G+w$.
\end{lem}


Let us consider a $C^{-1}$ biinvariant distance $\gamma$ on
$\mathcal{H}_c(\RN)$ which is invariant under the action of affine
at infinity Hamiltonians (such a condition is clearly satisfied by
Hofer's distance). For a sequence of Hamiltonian functions that are
affine at infinity with the same affine part, we can speak of its
limit for $\gamma$ by setting:
$$(\phi_{H_n+u})\stackrel{\gamma}{\to}\phi_{H+u} \text{ if and only if
}\gamma((\phi_{H+u})^{-1}\phi_{H_n+u},Id)\to 0.$$
Moreover, if
$(\phi_{H_n+u})\stackrel{\gamma}{\to}\phi_{H+u}$ and
$(\phi_{K_n+v})\stackrel{\gamma}{\to}\phi_{K+v}$ then
$$(\phi_{H_n+u}\phi_{K_n+v})\stackrel{\gamma}{\to}\phi_{H+u}\phi_{K+v}.$$
Indeed,  we have
\begin{eqnarray*}\lefteqn{\gamma((\phi_{H_n+u}\phi_{K_n+v})^{-1}(\phi_{H+u}\phi_{K+v}),Id)}\\
 & = &
 \gamma(\phi_{K+v}^{-1}(\phi_{H+u}^{-1}\phi_{H_n+u})\phi_{K+v}(\phi_{K+v}^{-1}\phi_{K_n+v}),Id)\\
 & \leq &
 \gamma(\phi_{H+u}^{-1}\phi_{H_n+u},Id)+\gamma(\phi_{K+v}^{-1}\phi_{K_n+v},Id).
\end{eqnarray*}
Finally notice that if $\|H_n-H\|_{C^0}\to 0$, then
$\phi_{H_n+u}\stackrel{\gamma}{\to} \phi_{H+u}$.

We are now ready for our proofs.

\medskip
\noindent\textit{Proof of lemma \ref{lineaire1}.} We just adapt the
proof of Cardin and Viterbo \cite{CV} to the "affine at infinity"
case.

First remark that the assumptions imply $\{u,v\}=0$. Then, a simple
computation shows that the flow
$$\psi_n^t=\phi_{H_n+u}^t\phi_{K_n+v}^s\phi_{H_n+u}^{-t}\phi_{K_n+v}^{-s}$$
is generated by the Hamiltonian function affine at infinity
$$\int_0^s\{H_n+u,K_n+v\}(\phi_{K_n+v}^{\sigma}\phi_{H_n+u}^{t}(x))d\sigma,$$
which $C^0$-converges to $0=\{u,v\}$ by assumption. Therefore,
$\psi_n^t$ converges for any $s$ and any $t$ to $Id$. But on the
another hand, according to the above remark, it converges to
$\phi_{H+u}^t\phi_{K+v}^s\phi_{H+u}^{-t}\phi_{K+v}^{-s}$. Hence
$\phi_{H+u}^t\phi_{K+v}^s\phi_{H+u}^{-t}\phi_{K+v}^{-s}=Id$ which
proves $\{H+u,K+v\}=0$.$\quad\Box$

\medskip
\noindent\textit{Proof of lemma \ref{lineaire2}.} First notice that
the assumptions imply $\{u,v\}=w$, $\{u,w\}=0$ and $\{v,w\}=0$, and
that the equalities $\{H+u,G+w\}=0$, $\{K+v,G+w\}=0$ follow from
lemma \ref{lineaire1}. Here we consider the flow
$$\psi_n^t=\phi_{G_n+w}^{-ts}\phi_{H_n+u}^t\phi_{K_n+v}^s\phi_{H_n+u}^{-t}\phi_{K_n+v}^{-s}$$
which is generated by
$$\left(-s(G_n+w)+\int_0^s\{H_n+u,K_n+v\}(\phi_{K_n+v}^{\sigma}\phi_{H_n+u}^{t})d\sigma\right)\circ\phi_{G_n+w}^{ts}.$$
This expression can be written
$$\left(\int_0^s(A_n+B_n)d\sigma\right)\circ\phi_{G_n+w}^{ts},$$
where $A_n=G_n-G_n(\phi_{K_n+v}^{\sigma}\phi_{H_n+u}^{t})$ and
$B_n=(\{H_n+u,K_n+v\}-(G_n+w))(\phi_{K_n+v}^{\sigma}\phi_{H_n+u}^{t})$.

By assumption, $B_n$ $C^0$-converges to $0$ and $A_n$ can be
written:
\begin{eqnarray*}
A_n & = & (G_n-G_n(\phi_{H_n+u}^{t})) + (G_n-G_n(\phi_{K_n+v}^{\sigma}))\circ\phi_{H_n+u}^{t}\\
 & = & \int_0^{t}\{G_n, H_n+u\}d\tau+\left(\int_0^{\sigma}\{G_n,
 K_n+v\}d\tau\right)\circ\phi_{H_n+u}^{t}\\
 & = & \int_0^{t}\{G_n+w, H_n+u\}d\tau+\left(\int_0^{\sigma}\{G_n+w,
 K_n+v\}d\tau\right)\circ\phi_{H_n+u}^{t},
\end{eqnarray*}
which implies that $A_n$ $C^0$-converges to $0$ too. It follows that
the generating Hamiltonian of $\psi_n^t$ $C^0$-converges to $0$, and
hence that $\psi_n^t$ $\gamma$-converges to $Id$. Since it also
converges to
$\psi^t:=\phi_{G+w}^{-ts}\phi_{H+u}^t\phi_{K+v}^s\phi_{H+u}^{-t}\phi_{K+v}^{-s}$,
we get $\psi^t=Id$ for any $s$ and $t$. Thus, the generating
Hamiltonian of $\psi_t$ vanishes identically:
$$\left(-s(G+w)+\int_0^s\{H+u,K+v\}(\phi_{K+v}^{\sigma}\phi_{H+u}^{t})d\sigma\right)\circ\phi_{G+w}^{ts}=0.$$
But since $G+w$ commutes with $H+U$ and $K+v$, we get:
$$\int_0^s(\{H+u,K+v\}-(G+w))(\phi_{K+v}^{\sigma}\phi_{H+u}^{t})d\sigma=0.$$
Taking derivative with respect to $s$, we obtain
$\{H+u,K+v\}-(G+w)=0$.$\quad\Box$

\section{Few additional remarks using the theory of distributions.}

The following results on Poisson brackets are obtained with the help
of distributions. No assumptions are made on the Lie algebra
generated by the Hamiltonian functions. They show in a certain way
why it is difficult to find examples of pseudo-representations whose
limit is not a representation.

\begin{prop}\label{distributions1} If $F_n$ $C^2$-converges to $F$ and $G_n$
$C^0$-converges to $G$. Then, $\{F_n,G_n\}$  converges to $\{F,G\}$
in the sense of distributions. As a consequence, if $\{F_n,G_n\}$
$C^0$-converges to $H$, then $\{F,G\}=H$.
\end{prop}

\medskip
\noindent\textit{Proof. } For any smooth compactly supported
function $\phi$,
\begin{eqnarray*}
\langle \{F_n,G_n\},\phi\rangle & =& \int
\derive{G_n}{q}\derive{F_n}{p}\phi-\int
\derive{G_n}{p}\derive{F_n}{q}\phi\\
& = & - \int G_n \derive{}{q}\left(\derive{F_n}{p}\phi\right)+\int
G_n\derive{}{p}\left(\derive{F_n}{q}\phi\right).
\end{eqnarray*}
By assumption, the integrands $C^0$-converge and hence the integrals
converge to $- \int G
\derive{}{q}\left(\derive{F}{p}\phi\right)+\int
G\derive{}{p}\left(\derive{F}{q}\phi\right)$ which equals
$\langle\{F,G\},\phi\rangle$.$\quad\Box$

\begin{prop}\label{distributions2}If $F_n$ $C^0$-converges to $F$, $G_n$
$C^0$-converges to $G$ and $\{F_p,G_q\}$  $C^0$-converges to $H$
when $p$ and $q$ go to infinity, then $\{F,G\}=H$.
\end{prop}

\medskip
\noindent\textit{Proof. } Take once again a compactly supported
smooth function $\phi$. Write
$$\langle\{F_p,G_q\}-\{F,G\},\phi\rangle=\langle\{F_p-F,G_q\},\phi\rangle+\langle\{F,G_q-G\},\phi\rangle.$$
By Proposition \ref{distributions1}, the first term converges to 0.
Hence for all $\eps>0$, there exists an integer $q_0$ such that for
any $q>q_0$, $|\langle\{F,G_q-G\},\phi\rangle|\leq\eps$.

Similarly, for each fixed $q$, there exists an integer $p_0$ such
that for any $p>p_0$, $|\langle\{F_p-F,G_q\},\phi\rangle|\leq\eps$.

Therefore, for all $\eps$ and all integers $p_1, q_1$, we can find
$p>p_1, q>q_1$ such that
$|\langle\{F_p,G_q\}-\{F,G\},\phi\rangle|\leq 2\eps$.

Thus we can construct two extractions $\chi,\psi$ such that
$\langle\{F_{\chi(n)},G_{\psi(n)}\}-\{F,G\},\phi\rangle$ converges
to 0. Since we have
$\langle\{F_{\chi(n)},G_{\psi(n)}\}-H,\phi\rangle\to 0$, it implies
$\langle\{F,G\},\phi\rangle=\langle H,\phi\rangle$, and this
equality holds for any $\phi$.$\quad\Box$

\bigskip

\noindent \textbf{\large{Acknowledgments. }} I warmly thank my
supervisor Claude Viterbo for all his advices and for hours of
fruitful discussion. I also thank Nicolas Roy for innumerable
interesting conversations on multiple subjects.

\nocite{*}
\bibliographystyle{plain}
\bibliography{article2}

\begin{thebibliography}{10}

\bibitem{CV}
\textsc{Cardin} F. and \textsc{Viterbo} C.
\newblock \textit{Commuting {H}amiltonians and {H}amilton-{J}acobi multi-time
  equations}.
\newblock preprint, math.SG/0507418.

\bibitem{EPZ}
\textsc{Entov} M.~\textsc{Polterovich} L. and \textsc{Zapolsky} F.
\newblock \textit{Quasi-morphisms and the Poisson bracket}.
\newblock preprint, math.SG/0605406, 2006.

\bibitem{FHV}
\textsc{Hofer}~H. \textsc{Floer} A. and \textsc{Viterbo} C.
\newblock The {W}einstein conjecture in ${P}\times$\textbf{C}$^l$.
\newblock {\em Math Z.}, (203):469--482, 1990.

\bibitem{H1}
\textsc{Hofer} H.
\newblock On the topological properties of symplectic maps.
\newblock {\em Proc. Roy. Soc. Edinburgh Sect. A}, 115:25--38, 1990.

\bibitem{HV1}
\textsc{Hofer} H. and \textsc{Viterbo} C.
\newblock The {W}einstein conjecture in cotangent bundles and related results.
\newblock {\em Annali Sc. Norm. Sup. Pisa}, 15, 1988.

\bibitem{HV2}
\textsc{Hofer} H. and \textsc{Viterbo} C.
\newblock The {W}einstein conjecture in the presence of holomorphic spheres.
\newblock {\em Comm. Pure Appl. Math.}, 1992.

\bibitem{HZ}
\textsc{Hofer} H. and \textsc{Zehnder} E.
\newblock {\em Symplectic invariants and Hamiltonian dynamics.}
\newblock Birkhauser, 1994.

\bibitem{Hu}
\textsc{Humilière} V.
\newblock \textit{On some completions of the space of {H}amiltonian maps}.
\newblock preprint, math.SG/0511418.

\bibitem{Kl}
\textsc{Klingenberg} W.
\newblock {\em Closed geodesics on Riemannian manifolds.}
\newblock Number~53 in CBMS Regional Conference Series in Mathematics.
  Published for the Conference Board of the Mathematical Sciences, Washington,
  DC; by the American Mathematical Society, Providence, RI, 1983.

\bibitem{LT}
\textsc{Liu} G. and \textsc{Tian} G.
\newblock Weinstein conjecture and {G}{W} invariants.
\newblock {\em Commun. Contemp. Math.}, 2, 2000.

\bibitem{Lu}
\textsc{Lu} G.
\newblock The weinstein conjecture in the uniruled manifolds.
\newblock {\em Math. Res. Lett.}, 7, 2000.

\bibitem{Sc}
\textsc{Schwarz} M.
\newblock On the action spectrum for closed symplectically aspherical
  manifolds.
\newblock {\em Pacific J. Math.}, 193:419--461, 2000.

\bibitem{T}
\textsc{Taubes} C.~H.
\newblock \textit{The Seiberg-Witten equations and the Weinstein conjecture}.
\newblock preprint, math.SG/0611007, 2006.

\bibitem{V2}
\textsc{Viterbo} C.
\newblock A proof of {W}einstein conjecture in \textbf{R}$^{2n}$.
\newblock {\em Ann. Inst. Poincar{é}, Anal. Non Lin´eaire,}, 4, 1987.

\bibitem{V1}
\textsc{Viterbo} C.
\newblock Symplectic topology as the geometry of generating functions.
\newblock {\em Math. Annalen}, 292:685--710, 1992.

\bibitem{Z}
\textsc{Zapolsky} F.
\newblock \textit{Quasi-states and the Poisson bracket on surfaces}.
\newblock preprint, math.SG/0703121, 2007.

\end{thebibliography}

\end{document}